\documentclass[12pt,a4paper,draft]{article}

\usepackage[reqno]{amsmath}
\usepackage{amssymb,amsthm,upref,amscd}
\usepackage[T1]{fontenc}
\usepackage{times}
\usepackage{color}

\newtheorem{theorem}{Theorem}[section]

\newtheorem{lemma}[theorem]{Lemma}

\theoremstyle{definition}
\theoremstyle{remark}
\numberwithin{equation}{section}

\begin{document}

\title{\textbf{Comment On the Uniqueness of the Ground State Solutions of a Fractional NLS with a Harmonic Potential    \ \\
		}}
		
		\author{Hichem Hajaiej$^{\mathrm{a}}$\thanks{%
				Email: hhajaie@calstatela.edu}\  ,	
			Linjie Song$^{\mathrm{b,c,}}$\thanks{%
				The author is supported by CEMS. Email: songlinjie18@mails.ucas.edu.cn.}\  \\
			{\small $^{\mathrm{a}}$ Department of Mathematics, California State University at Los Angeles,}\\
			{\small Los Angeles, CA 90032, USA}\\
			{\small $^{\mathrm{b}}$Institute of Mathematics, AMSS, Academia Sinica,
				Beijing 100190, China}\\
			{\small $^{\mathrm{c}}$University of Chinese Academy of Science,
				Beijing 100049, China}
		}
		\date{}
		\maketitle

\begin{abstract}
  The uniqueness of the positive ground state solutions of fractional Shrodinger equations with a harmonic potential has not been covered by the breakthrough method developed in \cite{Lenzmann1, Lenzmann2}. It has remained an open question for years. \cite{HD} and \cite{SS} were quite recently able to prove the existence of the ground state solutions and to derive important properties but they failed to address the uniqueness. \cite{HD} left it as an open question, and \cite{SS}  run numerical simulations that were in favor of the uniqueness. Very recently, the authors of this note developed a general and unified method to prove the uniqueness of the ground state solutions of a large class of variational problems, they also exhibited many examples to which their approach applies. After the publication of \cite{HS}, some colleagues reached out to us to know whether our method applies to the fractional Shrodinger equation with a harmonic potential. In this short note, we will explain how it applies, and we will also state some existence/ non-existence and multiplicity solutions of the normalized solutions.
\end{abstract}

\section{Introduction}

 In our recent paper, \cite{HS}, we developed a general and abstract method to prove the uniqueness of
positive ground state solutions. We did exhibit many examples to which our approach applies. We also
pointed out in the Introduction that our techniques and ideas apply to many other operators, domains,
and types of nonlinearity. Among the numerous examples, our approach applies to the equation
\begin{equation} \label{eq1.1}
\left\{
\begin{array}{cc}
(-\Delta)^{s}u + V(|x|)u = \lambda u + |u|^{q-2}u \ in \ \mathbb{R}^{N}, \\
u(x) \rightarrow 0 \ as \ |x| \rightarrow +\infty,
\end{array}
\right.
\end{equation}
with $0 < s < 1$, $2 < q < 2_{s}^{\ast}$ and when the potential $V$ is Harmonic, $V(x) = |x|^2$. Very recently, \cite{Gou} also tackled this problem. When $N\geq1,$ $0<s<1,$ $\lambda<\lambda_{1,s},$ where $\lambda_{1,s}=\inf\sigma\{(-\Delta)^s+|x|^2\},$ and $2<p<2_s^*,$ he was able to show the uniqueness of positive ground state solutions for \eqref{eq1.1} with Harmonic potential. His main result can be obtained by following the discussions in the proof of \cite[Theorem 2]{SS}). Indeed, we can find a suitable path parameterized by $s$ connecting (\ref{eq1.1}) to the case when $s = 1$ (consider $K(u,s) = (-\Delta)^{s}u + |x|^{2}u - \lambda u - |u|^{q-2}u$) and then similar discussions to \cite[Theorem 4.28, Lemmas 4.29, 4.30]{HS} yield to uniqueness of positive ground state. This was exactly what he did to show this property in \cite[Theorem 1.5]{HS}

Furthermore, along the lines developed in \cite{HS}, we have a second approach to show the uniqueness. More precisely, that we can solve the uniqueness problem for general trapping potential under the conditions given by \cite{HD, SS} even when we do not know the uniqueness of ground state for $s = 1$. We are preparing an abstract framework, and the general trapping potentials (including harmonic potential) can be viewed as an example.

Also, we can apply our abstract framework in \cite{HS} to (\ref{eq1.1}) when $V(|x|) = |x|^{2}$ (or any other general trapping potentials) to show the existence, non-existence and multiplicity of the normalized solutions (see Section 2).

For the convenience of the reader, we collect these results below.

\section{Existence/Non-existene/Multiplicity of normalized solutions}

Consider the following equation: (\eqref{eq1.1} with $V(x)=|x|^2$)
\begin{equation} \label{eq2.1}
\left\{
\begin{array}{cc}
(-\Delta)^{s}u + |x|^{2}u = \lambda u + |u|^{q-2}u \ in \ \mathbb{R}^{N}, \\
u(x) \rightarrow 0 \ as \ |x| \rightarrow +\infty,
\end{array}
\right.
\end{equation}
with $0 < s \leq 1$, $2 + 4s/N < q < 2_{s}^{\ast}$.

Apply the abstract framework in \cite{HS} with
$$
W = \{u \in H^{s}(\mathbb{R}^{N}): \int_{\mathbb{R}^{N}}|x|^{2}u^{2} < +\infty\},
$$
$$
\Phi_{\lambda}(u) = \frac{1}{2}\int_{\mathbb{R}^{N}}[|(-\Delta)^{\frac{s}{2}}u|^{2} + |x|^{2}u^{2} - \lambda u^{2}]dx - \frac{1}{q}\int_{\mathbb{R}^{N}}|u|^{q}dx.
$$

\begin{lemma}
(Pohoz$\tilde{a}$ev identity) Let $u \in W$ be a weak solution of (\ref{eq2.1}), then the following integral identity holds true
\begin{eqnarray} \label{eq2.2}
	(N-2s)\int_{\mathbb{R}^{N}}|(-\Delta)^{\frac{s}{2}} u|^{2}dx + (N+2)\int_{\mathbb{R}^{N}}|x|^{2}u^{2}dx = N\lambda\int_{\mathbb{R}^{N}}u^{2}dx + \frac{2N}{q}\int_{\mathbb{R}^{N}}|u|^{q}dx.
\end{eqnarray}
\end{lemma}

Let $u_{\lambda}$ be the unique ground state of (\ref{eq2.1}) for $\lambda < \lambda_{1}$ where $\lambda_{1} = \inf \sigma((-\Delta)^{s} + |x|^{2})$. Similar to \cite[Lemma 4.22]{HS}, we know that $(\lambda,u_{\lambda})$ is a $C^{1}$ curve in $\mathbb{R} \times W$. Then, we analyze the behavior of $\int_{\mathbb{R}^{N}}u^{2}dx$ when $\lambda \rightarrow -\infty$ and $\lambda \rightarrow 0$.

\begin{lemma}
Let $u_{\lambda}$ be the unique ground state solution of (\ref{eq2.1}). Assume that $2 + 4s/N < q < 2_{s}^{\ast}$. Then $\int_{\mathbb{R}^{N}}u_{\lambda}^{2}dx \rightarrow 0$ as $\lambda \rightarrow -\infty$.
\end{lemma}

\textit{Proof.  } Since $u_{\lambda}$ solves (\ref{eq2.1}), we have
$$
\int_{\mathbb{R}^{N}}|(-\Delta)^{\frac{s}{2}} u_{\lambda}|^{2}dx + \int_{\mathbb{R}^{N}}|x|^{2}u_{\lambda}^{2}dx = \lambda\int_{\mathbb{R}^{N}}u_{\lambda}^{2}dx + \int_{\mathbb{R}^{N}}|u_{\lambda}|^{q}dx.
$$
The latter, together with (\ref{eq2.2}), imply that
$$
s\int_{\mathbb{R}^{N}}|(-\Delta)^{\frac{s}{2}} u_{\lambda}|^{2}dx = \int_{\mathbb{R}^{N}}|x|^{2}u_{\lambda}^{2}dx + \frac{q-2}{2q}N\int_{\mathbb{R}^{N}}|u_{\lambda}|^{q}dx,
$$
$$
2(1+s)\int_{\mathbb{R}^{N}}|x|^{2}u_{\lambda}^{2}dx - 2s\lambda\int_{\mathbb{R}^{N}}u_{\lambda}^{2}dx = [\frac{2N}{q}-(N-2s)]\int_{\mathbb{R}^{N}}|u_{\lambda}|^{q}dx
$$
yielding to
\begin{eqnarray}
	\Phi_{\lambda}(u_{\lambda}) &=& \frac{1+s}{2s}(1 + k)\int_{\mathbb{R}^{N}}|x|^{2}u_{\lambda}^{2}dx - (1 + k)\frac{\lambda}{2}\int_{\mathbb{R}^{N}}u_{\lambda}^{2}dx \nonumber \\
	&\geq&  - (1 + k)\frac{\lambda}{2}\int_{\mathbb{R}^{N}}u_{\lambda}^{2}dx,
\end{eqnarray}
where $k = \frac{(q-2)N-4s}{2N-(N-2s)q}$. Since $2 + 4s/N < q < 2_{s}^{\ast}$, we have that $k > 0$. Then as a corollary of \cite[Theorem 2.2 (ii)]{HS}, we can deduce that $\int_{\mathbb{R}^{N}}u_{\lambda}^{2}dx \rightarrow 0$ as $\lambda \rightarrow -\infty$.
\qed\vskip 5pt

\begin{theorem}
	Assume that $2 + 4s/N < q < 2_{s}^{\ast}$. Then there exists $c_{0}$ such that
	
	$(i)$ if $c < c_{0}$, (\ref{eq2.1}) has at least two ground states $u_{\lambda}$, $u_{\widetilde{\lambda}}$ with $\lambda < \widetilde{\lambda} < \lambda_{1}$, $\int_{\mathbb{R}^{N}}u_{\lambda}^{2}dx = \int_{\mathbb{R}^{N}}u_{\widetilde{\lambda}}^{2}dx = c$, and for almost every $c \in (0,c_{0})$, $u_{\lambda}$ is orbitally unstable while $u_{\widetilde{\lambda}}$ is orbitally stable;
	
	$(ii)$ (\ref{eq2.1}) has at least one ground state $u_{\lambda}$ with $\lambda < \lambda_{1}$, $\int_{\mathbb{R}^{N}}u_{\lambda}^{2}dx = c_{0}$;
	
	$(iii)$ if $c > c_{0}$, (\ref{eq2.1}) has no ground state $u_{\lambda}$ with  $\int_{\mathbb{R}^{N}}u_{\lambda}^{2}dx = c$.
\end{theorem}

\textit{Proof.  } $\lambda_{1}$ is the first eigenvalue of $(-\Delta)^{s} + |x|^{2}$ and simple. By standard bifurcation arguments, we can show that $u_{\lambda} \rightarrow 0$ in $W$ as $\lambda \rightarrow \lambda_{1}$ and thus $\int_{\mathbb{R}^{N}}u_{\lambda}^{2}dx \rightarrow 0$. The rest of the proof can be completed in a similar way to \cite[Theorema 4.1, 4.21]{HS}.
\qed


\vskip 5pt

\begin{thebibliography}{99}

\bibitem{Lenzmann1}   R. L. Frank, E. Lenzmann, Uniqueness of non-linear ground states for fractional Laplacians in $\mathbb{R}$,
R. Acta Math. 210 (2013), no. 2, 261–318.

\bibitem{Lenzmann2}  R. L. Frank, E. Lenzmann, L. Silvestre, Uniqueness of radial solutions for the fractional Laplacian, Comm. Pure Appl. Math. 69 (2016), no. 9, 1671–1726.

\bibitem{HD} H. Hajaiej, Z Ding,   On a fractional Schr\"{o}dinger equation in the presence of Harmonic potential, Electronic Research Archive 2021, doi: 10.3934/era.2021047


\bibitem{HS}  H. Hajaiej, L.J. Song, A general and unified method to prove the existence of normalized solutions and some applications, arXiv.2208.11862.

\bibitem{SS}  M. Stanislavova, A.G. Stefanov, Ground states for the nonlinear Schr\"{o}dinger equation under a general trapping potential, J. Evol. Equ. 21 (2021), no. 1, 671-697

\end{thebibliography}
\end{document}